\documentclass[12pt]{amsart}
\newif\ifpdf\ifx\pdfoutput\undefined\pdffalse\else\pdfoutput=1\pdftrue\fi

\usepackage{hyperref}
\usepackage{graphicx,amsmath,amsfonts,latexsym,amssymb,amsthm,amscd}
\usepackage{latexsym,amscd}

  \newcommand{\field}[1]{\mathbb{#1}}
  \newcommand{\rbb}{\field{R}}
  \newcommand{\cbb}{\field{C}}

  \newcommand{\Spin}{\mathrm{Spin}}
  \newcommand{\SO}{\mathrm{SO}}
  \newcommand{\I}{\,\mathrm{i}\,}
  \newcommand{\PsDO}{\Psi\mathrm{DO}}

  \newcommand{\tr}{\mathrm{Tr}}

  \newcommand{\beq}{\begin{equation}}
  \newcommand{\eeq}{\end{equation}}

\newcommand{\reals}{{\bf R}}

  \newtheorem{theorem}{Theorem}[section]
  
  \newtheorem{lem}[theorem]{Lemma}
  \newtheorem{pro}[theorem]{Proposition}
  \newtheorem{cor}[theorem]{Corollary}

\begin{document}

\title[High energy limits and frame flows]{High energy limits of Laplace-type and
Dirac-type eigenfunctions and frame flows}

\author[D. Jakobson]{Dmitry Jakobson}
\address{Department of Mathematics and Statistics, McGill
University, 805 Sherbrooke Str. West, Montr\'eal QC H3A 2K6,
Ca\-na\-da.} \email{jakobson@math.mcgill.ca}

\author[A. Strohmaier]{Alexander Strohmaier}
\address{Mathematisches Institut, Universit\"at Bonn, Beringstrasse 1,
D-53115 Bonn, Germany} \email{strohmai@math.uni-bonn.de}



\date{\today}

\thanks{The first author was supported by NSERC, FQRNT and Dawson fellowship.}

\maketitle

\begin{abstract}
 We relate high-energy limits of Laplace-type and Dirac-type operators
 to frame flows on the corresponding manifolds, and show that the
 ergodicity of frame flows implies quantum ergodicity in an appropriate sense
 for those operators. Observables for the corresponding quantum systems
 are matrix-valued pseudodifferential operators and therefore
 the system remains non-commutative in the high-energy limit. We discuss
 to what extent the space of stationary high-energy states behaves classically.
\end{abstract}

\section{Introduction and main results}

If $X$ is an oriented closed Riemannian manifold and $\Delta$ the Laplace operator
on $X$, then a complete orthonormal sequence of eigenfunctions $\phi_j \in
L^2(X)$ with eigenvalues $\lambda_j \nearrow \infty$ is known to converge
in the mean to the Liouville measure, in the sense that
\begin{gather*}
 \lim_{N \to \infty} \frac{1}{N} \sum_{j \leq N} \langle \phi_j,
 A \phi_j \rangle = \int_{T_1^*X} \sigma_A(\xi) dL(\xi),
\end{gather*}
for any zero order pseudodifferential operator $A$,
where integration is with respect to the normalized Liouville measure on the
unit cotangent bundle $T_1^* X$, and $\sigma_A$ is the principal symbol of $A$.
In particular, $A$ might be a smooth function on $X$ and the above implies that
the sequence
\begin{gather*}
 \frac{1}{N} \sum_{j \leq N} |\phi_j(x)|^2
\end{gather*}
converges to the normalized Riemannian measure in the weak topology of measures.
In case the geodesic flow on $T_1 X$ is ergodic it is known that the following stronger
result holds.
\begin{gather*}
 \lim_{N \to \infty }\frac{1}{N} \sum_{j \leq N} |\langle \phi_j, A \phi_j \rangle -
 \int_{T_1^* X} \sigma_A(\xi) dL(\xi)| =0.
\end{gather*}
This property is commonly referred to as quantum ergodicity and it is
equivalent
to the existence of a density-one-subsequence $\phi_j'$ such that
\begin{gather*}
 \lim_{j \to \infty} \langle \phi_j', A \phi_j' \rangle = \int_{T_1^* X}
\sigma_A(\xi) dL(\xi),
\end{gather*}
for any zero order pseudodifferential operator $A$ (see
\cite{MR0402834, MR1239173, MR0818831, MR0916129}).

We show in this paper that the high energy behavior of the Dirac
operator $D$ acting on spinors on a closed spin manifold $X$ is
determined by the frame flow in the same manner, as the geodesic
flow determines the high energy limit of the Laplace operator. If
$F_k X$ is the bundle of oriented orthonormal $k$-frames in $T^* X$,
then projection to the first vector makes $F_k X \to T_1^* X$ into a
fiber bundle. In particular for $k=n$ this is the full frame bundle
and $FX=F_n X$ is a principal fiber bundle over $T_1^*X$ with
structure group $SO(n-1)$. Transporting covectors parallel along
geodesics extends the Hamiltonian flow on $T_1^*X$ to a flow on
$F_kX$. This is the so-called $k$-frame flow. In case $k=n$ we will
refer to it simply as the frame flow. Of course ergodicity of the
$k$-frame flow for any $k$ implies ergodicity of the geodesic flow,
whereas the conclusion in the other direction is not always true
(cf. section \ref{sec:examples}). Still there are many examples investigated in the
literature when the frame flow is ergodic. Our first main result is,
that quantum ergodicity holds for eigensections of the Dirac
operator in case the frame flow is ergodic.
\begin{theorem}\label{erg:dirac}
 Let $X$ be a closed Riemannian spin manifold of dimension $n \geq 3$
 with Dirac operator $D$ acting on sections of the spinor bundle.
 Suppose that $\phi_j \in L^2(X;S)$ is an orthonormal sequence of eigensections
 of $D$ with eigenvalues $\lambda_k \nearrow \infty$ such that
 the $\phi_k$ span\footnote{in the sense that the linear hull is dense}
 the positive energy subspace of $D$. Then, if the frame flow on $FX$
 is ergodic, we have
 \begin{gather*}
  \lim_{N \to \infty }\frac{1}{N} \sum_{j \leq N} |\langle \phi_j, A \phi_j \rangle -
 \frac{1}{2^{[\frac{n}{2}]}}\int_{T_1^* X}
 \tr{\left(\left(1+\gamma(\xi)\right)\sigma_A(\xi) \right)} dL(\xi)| =0,
 \end{gather*}
 for all $A \in \PsDO_{cl}^0(X,S)$. Here $\gamma(\xi)$ denotes the operator of Clifford
 multiplication with $\xi$.
 In particular there is a density one subsequence $\phi_j'$ such that
 \begin{gather*}
  \langle \phi_j', A \phi_j' \rangle \to \frac{1}{2^{[\frac{n}{2}]}}\int_{T_1^* X}
  \tr{\left(\left(1+\gamma(\xi)\right) \sigma_A(\xi)\right)} dL(\xi).
 \end{gather*}
 A similar statement holds for the negative energy subspace.
\end{theorem}

Another result is that the $\left(2\,
\mathrm{min}(p,n-p)\right)$-frame  flow determines the high energy
behavior of the Laplace-Beltrami operator $\Delta_p$ acting on the space
$C^\infty(X,\Lambda^{p}_\cbb X)$ of
complex-valued $p$-forms. Note that the Hodge decomposition implies that there are
three invariant subspaces for $\Delta_p$, namely the closures of $d
C^\infty(X,\Lambda^{p-1}_\cbb X)$, $\delta C^\infty(X,\Lambda^{p+1}_\cbb X)$
and the finite dimensional space of harmonic forms
The latter
subspace plays no role for the high energy behavior. The eigenspaces
of the first subspace consist of exterior derivatives of $p-1$
eigenforms their high energy behavior is therefore determined by the
high energy behavior of $\Delta_{p-1}$. In particular the high
energy behavior of $\Delta_1$ restricted to $d C^\infty(X)$ is
controlled by the geodesic flow. We therefore look at the second
subspace only. Note that any coclosed form which is a nonzero
eigenvalue to $\Delta_p$ is coexact. Hence, to investigate the high
energy behavior of $\Delta_p$ we have to look at the system
\begin{gather*}
 \Delta_p \phi_j = \lambda_j \phi_j,\\
 \delta \phi_j = 0.
\end{gather*}
In case $p=1$ such systems appear in physics if one investigates electromagnetic fields
(Maxwell's equations) or the Proca equation for spin 1 particles. The
restriction to the coclosed forms corresponds to a gauge condition
which restricts to the transversal subspaces. Our result is,
that this system is quantum ergodic, if the $2\,
\mathrm{min}(p,n-p)$-frame flow is ergodic.

\begin{theorem}\label{erg:forms}
 Let $X$ be an oriented closed Riemannian manifold of dimension $n \geq 3$
 and let $0 < p < n$. Suppose that $\phi_k$ is an
 orthonormal sequence of eigen-$p$-forms satisfying
 \begin{gather*}
 \Delta_p \phi_k = \lambda_k \phi_k,\\
 \delta \phi_k =0,
\end{gather*}
such that the $\phi_k$ span $\mathrm{ker} (\delta)$ and with $\lambda_k \nearrow \infty$.
Suppose that $p \not= \frac{n-1}{2}$. Then, if the
$\left(2\, \mathrm{min}(p,n-p)\right)$-frame flow is ergodic,
the system is quantum ergodic in the sense that
\begin{gather*}
  \lim_{N \to \infty }\frac{1}{N} \sum_{k \leq N} |\langle \phi_k, A \phi_k \rangle -
 \omega_t(\sigma_A)| =0,
 \end{gather*}
 for all $A \in \PsDO_{cl}^0(X;\Lambda^p_\cbb X)$. In particular there
is a density one subsequence
 $\phi_k'$ such that
 \begin{gather*}
  \lim_{k \to \infty} \langle \phi_k', A \phi_k' \rangle
  = \omega_t(\sigma_A), \quad \textrm{for all}\; A \in \PsDO_{cl}^0(X;\Lambda^p_\cbb X).
 \end{gather*}
 Here $\omega_t$ is a state on the $C^*$-algebra of continuous
 $\mathrm{End}(\Lambda^p_\cbb X)$-valued functions
 on $T_1^*X$ which is defined by
\begin{gather*}
 \omega_t(a):={n-1 \choose p}^{-1} \int_{T_1^*X} \tr \left(i(\xi)i^*(\xi)
a(\xi)\right) dL(\xi),
\end{gather*}
where $i(\xi)$ is the operator of interior multiplication with
$\xi$, and the adjoint $i^*(\xi)$ is the operator of exterior
multiplication with $\xi$.
\end{theorem}

\noindent
Note that the system
\begin{gather*}
 \Delta_k \phi_j = \lambda_j \phi_j,\\
 d \phi_j = 0.
\end{gather*}
is equivalent to our system with $p=n-k$ via the Hodge star operator.

The restriction $p \not= \frac{n-1}{2}$ is necessary since if
$p = \frac{n-1}{2}$ the operator
$\I^{p+1}\delta *$ leaves the space $\overline{\mathrm{Rg}(\delta)}$
invariant and commutes with $\Delta_p$ ($*$ is the Hodge star operator).
In this case our result is

\begin{theorem}\label{erg:forms2}
 Let $X$ be an oriented closed Riemannian manifold of odd dimension $n \geq 3$.
 Let $p=\frac{n-1}{2}$ and
 suppose that $\phi_k$ is an orthonormal sequence of eigen-$p$-forms satisfying
 \begin{gather*}
 \Delta_p \phi_k = \lambda_k \phi_k,\\
 \delta \phi_k =0,\\
 \I^{p+1} \delta * \phi_k = \pm \sqrt{\lambda_k} \phi_k
\end{gather*}
such that the $\phi_k$ span $\overline{\mathrm{Ran}(\delta \pm \I^{p+1}
\Delta_p^{-1/2}\delta * \delta)}$
and with $\lambda_k \nearrow \infty$.
Then, if the $(n-1)$-frame flow is ergodic,
the system is quantum ergodic in the sense that
\begin{gather*}
  \lim_{N \to \infty }\frac{1}{N} \sum_{k \leq N} |\langle \phi_k, A \phi_k \rangle -
 \omega_\pm(\sigma_A)| =0,
 \end{gather*}
 for all $A \in \PsDO_{cl}^0(X;\Lambda^p_\cbb X)$. In particular there is a
density one subsequence
 $\phi_k'$ such that
 \begin{gather*}
  \lim_{k \to \infty} \langle \phi_k', A \phi_k' \rangle
  = \omega_\pm(\sigma_A), \quad \textrm{for all}\; A \in \PsDO_{cl}^0(X;\Lambda^p_\cbb X).
 \end{gather*}
 Here the states $\omega_\pm$ are defined by
\begin{gather*}
 \omega_\pm(a):=\frac{(p!)^2}{(2p)!} \int_{T_1^*X} \tr
 \left( \left(1 \pm \I^{p} i(\xi) *\right) i(\xi) i^*(\xi) a(\xi)\right) dL(\xi).
\end{gather*}
\end{theorem}

As the example of K\"ahler manifolds (see section \ref{sec:examples}) 
show the above theorems do not hold if we assume ergodicity of the geodesic
flow only.

Our analysis is based on a version of Egorov's theorem for matrix valued
operators. A second order differential
operator $P$ acting on sections of a vector bundle $E$
is said to be of Laplace type if $\sigma_P(\xi)=g(\xi,\xi )\;\mathrm{id}_E$,
 i.e. if in local coordinates it is of the form $P=-\sum_{i,k} g^{ik}
\partial_i \partial_k + \; \textrm{lower order terms}$.
Examples are the Laplace-Beltrami operator $\Delta_p$ acting on
$p$-forms or the square $D^2$ of the Dirac operator on a Riemannian
spin manifold. For such operators the first order term (the
subprincipal symbol) defines a connection $\nabla^E$ on the bundle
$E$. We will prove a Egorov theorem for matrix-valued
pseudodifferential operators acting on sections of $E$. More
precisely, for $A \in \PsDO_{cl}^0(X,E)$, a zero order classical
pseudodifferential operator, the principal symbol $\sigma_A$ is an
element in $C^\infty(T^*_1 X,\mathrm{End}(\pi^*(E)))$, where
$\pi^*(E)$ is the pull-back of the bundle $E \to X$ under the
projection $\pi: T_1^*X \to X$. Note that the connection $\nabla^E$
determines a connection $\nabla$ on $\mathrm{End}(\pi^*(E))$.
Parallel transport along the Hamiltonian flow of $\sigma_P$ then
determines a flow $\beta_t$ acting on $C^\infty(T^*_1
X,\mathrm{End}(\pi^*(E)))$. Our version of Egorov's theorem
specialized to Laplace type operators reads as follows.
\begin{theorem}
 If $A \in \PsDO_{cl}^0(X,E)$ and if $P$ is a
 positive second order differential operator of Laplace type then for all $t \in \mathbb{R}$
 the operators
 $A_t:=e^{+\I t P^{1/2}} A e^{-\I t P^{1/2}}$ are again in $\PsDO_{cl}^0(X,E)$
 and $\sigma_{A_t}=\beta_t(\sigma_A)$.
\end{theorem}
We actually prove a more general version of this theorem which applies
to flows generated by first order pseudodifferential operators with
real scalar principal symbols.
Note that unlike in the scalar case
the first order terms are needed to determine the flow. We show that for
pseudodifferential operators with real scalar principal part the subprincipal
symbol is invariantly defined as a partial connection along the Hamiltonian
vector field, thus allowing us to define all flows without referring to local
coordinate systems.

\subsection{Discussion}
Dirac equation on $\reals^3$ (and, more generally, on $\reals^d$)
has been studied from the semiclassical point of view in the papers
\cite{MR1644088, MR1694732, MR1838911, MR2073612, BoG04.2} of Bolte,
Glaser and Keppeler.  The authors would like to thank J. Bolte for
bringing to their attention this problem on manifolds. Unlike in the
works of Bolte, Glaser and Keppeler we investigate the high energy
limit rather than the semiclassical limit. Therefore, all nontrivial
dynamical effects are due to the nontrivial curvature of the
spin-connection. This is conceptually different from the stated
previous results where quantum ergodicity is due to a spin
precession in an external magnetic field. External fields are not
seen in the high energy limit and therefore the strict analog of the
result of Shnirelman, Colin de Verdi\`ere and Zelditch
\cite{MR0402834, MR1239173, MR0818831, MR0916129} can not be
expected to hold in $\mathbb{R}^n$ or on manifolds with integrable
geodesic flow.

Apart from working on manifolds our methods also differ from
those previously employed as we take  the absolute value of the
Dirac operator instead of the Dirac operator
itself as the generator of the dynamics.
This has the advantage of allowing for the full algebra of matrix
valued functions as the observable algebra rather than a subalgebra.
One can also justify this from a physical point of view.
Namely, in a fully quantized theory the generator of the time evolution
on the $1$-particle Hilbert spaces is the absolute value of the Dirac
operator. Furthermore, on the electron $1$-particle subspace these
two operators coincide.

Our theorem \ref{erg:forms2} for $n=3$ and $k=1$ deals with the
electromagnetic field on a $3$-dimensional compact manifold. The
statement of theorem \ref{erg:forms2} means that quantum ergodicity
holds for circular polarized photons if the $2$-frame flow is
ergodic.

We would also like to mention that the Egorov theorem as we state it
is related to a work of Dencker (\cite{MR661876}),
who proved a propagation of singularity theorem for systems of real
principal type. It follows from his work that the polarization set of
solutions to the Dirac equation is invariant under a certain flow
similar to ours.
We also refer the reader to \cite{EW96}, \cite{GMMP97}
and references therein, and \cite{San99} for discussion of semiclassical
limits for matrix-valued operators, and relations to parallel transport.

Since the high energy limit of the Quantum system associated to Laplace
type operators on vector bundles is non-commutative the apropriate
language to investigate questions of ergodicity is the language of 
$C^*$-dynamical systems and states (see Appendix \ref{appa}). This was already advertized by S. Zelditch 
in \cite{MR1384146} and it is shown there that for a large class of
abstract $C^*$-dynamical systems classical ergodicity implies quantum
ergodicity. The assumptions under which the theorems are stated in
\cite{MR1384146} ($G$-abelianness or classical abelianness) are in general not
satisfied in the examples we study. The method of the proof can be adapted
to our situation, however.
In our work we identify the classical
flows corresponding to the Dirac operator and the Hodge Laplacian as
frame flows, which allows us to use the results obtained by Brin, Arnold,
Pesin, Gromov, Karcher, Burns and Pollicott to exhibit many examples
of manifolds where quantum ergodicity holds for Dirac operator, see Corollary
\ref{cor:examples}.  The connection to their work has not been made before in
the literature on quantum ergodicity.  Finally, our results on quantum
ergodicity for $p$-eigenforms for Hodge Laplacian, and the role played by
$2\min(p,n-p)$-frame flow seem to be completely new.  It is a hope
of the authors that their results will stimulate further studies of
relationship between ergodic theory of partially hyperbolic dynamical
systems, and high energy behavior of eigenfunctions of matrix-valued
operators.

\section{Ergodic frame flows: known examples}\label{sec:examples}
The $k$-frame flow $\Phi^t$, $k\geq 2$ is defined as follows: let
$(v_1,\ldots,v_k)$ be an ordered orthonormal set of $k$ unit vectors
in $T_pX$.  Then $\Phi^t v_1=G^t v_1$, where $G^t$ is the geodesic
flow.  $\Phi^t v_j, 2\leq j\leq k$ translates $v_j$ by the parallel
translation at distance $t$ along the geodesic determined by $v_1$.
Here we summarize the cases when the frame flow is known to be
ergodic. A $k$-frame flow is a ${\rm SO}(k-1)$-extension of the
geodesic flow; on an $n$-dimensional manifold, $k$-frame flow is a
factor of the $n$-frame flow for $2\leq k<n$, so ergodicity of the
latter implies ergodicity of the former.  Frame flow preserves
orientation, so in dimension $2$ its ergodicity (restricted to
positively-oriented frames, say) is equivalent to the ergodicity of
the geodesic flow.

Frame flows were considered by Arnold in \cite{MR0158330}.  In
negative curvature, they were studied by Brin, together with Gromov,
Karcher and Pesin, in a series of papers \cite{MR0343316, MR0370660,
MR0394764, MR0582702, MR0670078, MR0756723}. Recently, a lot of
progress was made in understanding ergodic behavior of general
partially hyperbolic systems, including frame flows.  In the current
paper, the authors are primarily interested in the ergodicity of the
flow; the most recent paper dealing with that question appears to be
\cite{MR1988429} by Burns and Pollicott, where the authors establish
ergodicity under certain pinched curvature assumptions in
``exceptional'' dimensions $7$ and $8$, see below.

In the sequel, we shall assume that $M$ is negatively curved with
sectional curvatures satisfying
$$
-K_2^2\leq K \leq -K_1^2.
$$
The frame flow is known to be ergodic and have the $K$ property
\begin{itemize}
\item[1)] if $M$ has constant curvature \cite{MR0394764, MR0343316};
\item[2)] for an open and dense set of negatively curved metrics (in
the $C^3$ topology) \cite{MR0370660};
\item[3)] if $n$ is odd, but not equal to $7$ \cite{MR0582702}; or if $n=7$ and
$K_1/K_2>0.99023...$ \cite{MR1988429};
\item[4)] if $n$ is even, but not equal to $8$, and $K_1/K_2>0.93$,
\cite{MR0756723}; or if $n=8$ and $K_1/K_2>0.99023...$
\cite{MR1988429}.
\end{itemize}

By Theorem \ref{erg:spin} and results of section \ref{sec:hodge}, we
have the following
\begin{cor}\label{cor:examples}
Quantum ergodicity for Dirac operator and for Hodge Laplacian
(conclusions of Theorems \ref{erg:dirac}, \ref{erg:forms} and
\ref{erg:forms2}) hold in each of the cases (1)-(4).
\end{cor}

The frame flow is {\em not} ergodic on negatively-curved K\"ahler
manifolds, since the almost complex structure $J$ is preserved. This
is the only known example in negative curvature when the geodesic
flow is ergodic, but the frame flow is not. In fact, given an
orthonormal $k$-frame $(v_1,\ldots,v_k)$, the functions $(v_i,J
v_j), 1\leq i,j\leq k$ are first integrals of the frame flow, and in
some cases it is possible to describe the ergodic components,
\cite{MR0670078, MR0343316, MR0582702}. Note that the conclusion of
Theorem \ref{erg:forms} is false in the K\"ahler case, because
the decomposition into $(p,q)$-forms is a decomposition into
invariant subspaces of the Laplace-Beltrami operator. The K\"ahler
case is interesting in its own right and will be discussed in a
forthcoming paper.

The frame flow is conjectured to be ergodic whenever the curvature
satisfies $-1<K<-1/4$, cf. \cite{MR0670078}.  That conjecture is
still open.

\section{Microlocal analysis for operators on vector bundles}

\subsection{The subprincipal symbol}

Let $X$ be a closed manifold.
Suppose that $P \in \PsDO_{cl}^m(X,\Lambda^{n/2}X)$. Then the principal
symbol $\sigma_P$ is well defined as a function on the cotangent bundle
$\dot T^* X=T^*X \backslash 0$ which is smooth and positively homogeneous of
degree $m$.
The subprincipal symbol $\mathrm{sub}(P)$ is defined in local
coordinates by
\begin{gather} \label{subprince}
 \mathrm{sub}(P):= p_{m-1}-\frac{1}{2 \I} \sum_{j} \frac{\partial^2
 p_m}{\partial x^j \partial \xi_j}
\end{gather}
where the functions $p_m$ are the terms homogeneous of degree $m$ in
the asymptotic expansion of the full symbol of $P$. Surprisingly,
the subprincipal symbol turns out to be well defined as a function
on $\dot T^*X$ (see \cite{Hormander:1972}, ch 5.2). However, the
situation changes if we have $P \in \PsDO_{cl}^m(X,\Lambda^{n/2}X
\otimes E)$ for some vector bundle $E$. In this case the principal
symbol $\sigma_P$ is in $C^\infty(\dot T^* X,\mathrm{End}(\pi^*(E)))$,
whereas the subprincipal symbol has a more complicated
transformation law under a change of a bundle chart. More explicitly
we have with $\phi(x,\xi)=(x,\xi)$ for some section $w \in
C^\infty(X;E)$ using local coordinates and a local trivialization
(see \cite{Hormander:1972}, Equ. 5.2.2)
\begin{gather}\label{trans}
 e^{-\I \phi} P(e^{\I \phi}w) = p(x,\phi'_x)w
 - \frac{1}{2 \I} \sum_{j} \frac{\partial^2 p}{\partial x^j \partial \xi_j} w\\
 + \sum_{j} \left( p^{(j)}_m(x,\phi'_x) \frac{1}{\I} \frac{\partial w}{\partial x^j}+ \frac{1}{2 \I}
 \frac{\partial p^{(j)}_m(x,\phi'_x)}{\partial x^j} w \right) \quad \mathrm{mod} \nonumber
 \quad S^{m-2},
\end{gather}
where $p^{(j)}_m=\frac{\partial p_m}{\partial \xi_j}$. Suppose now that
$p_m=\sigma_P$ is scalar and real, i.e. $\sigma_P(x,\xi)=h(x,\xi) \mathrm{id}_{E_x}$
for some $h \in C^\infty(\dot T^* X,\mathbb{R})$.
The Hamiltonian vector field associated with the principal symbol
$\sigma_P$ of $P$ is a vector field on $\dot T^* X$ and defined in local
coordinated by
\begin{gather}
 H_P=\sum_j \left( \frac{\partial \sigma_P}{\partial \xi_j}\frac{\partial}{\partial x^j}
 -\frac{\partial \sigma_P}{\partial x^j}\frac{\partial}{\partial \xi_j} \right).
\end{gather}
Now $p^{(j)}_m(x,\phi'_x)$ has a nice interpretation in terms of
$H_P$. Namely, if $(x,\xi) \in \dot T^*X$, then $p^{(j)}_m(x,\xi)$
is the push forward of $H_P(x,\xi)$ under the projection $\pi: \dot T^*X \to X$
expressed in local coordinates, i.e.
\begin{gather}
 \sum_{j} p^{(j)}(x,\xi) \frac{\partial }{\partial x^j} = \pi_*(H_P(x,\xi)).
\end{gather}
Therefore, in case $E$ is trivial the last sum in (\ref{trans})
is exactly the Lie Derivative $-\I \mathcal{L}_{v} w$ of half densities
along the vector field $v=\pi_*(H_P(x,\phi'_x))$, which is defined without
reference to the local coordinate system and depends on the function $\phi$ only.
Hence, if we fix a local trivialization of $E$ and change coordinates on the
base manifold only, $\sigma_P$ and $\mathrm{sub(P)}$ transform as functions
in $C^\infty(\dot T^* X,\mathrm{End}(\pi^*(E)))$.
Note that $\sum_{j} p^{(j)}_m(x,\phi'_x) \frac{1}{\I}\frac{\partial w}{\partial x^j}$ is the
only term in (\ref{trans}) which depends on the derivative of $w$.
Hence, under a change of bundle charts by the local function
$A \in C^\infty(X,\mathrm{GL}(E))$ we get the transformation law
\begin{gather}\label{transform}
 \mathrm{sub}(P) \to A^{-1} \mathrm{sub}(P) A + A^{-1} p^{(j)} \partial_j A=\\= \nonumber
 A^{-1} \mathrm{sub}(P) A + A^{-1} \frac{1}{\I}H_P A
\end{gather}
whereas $\sigma_P$ transforms as a function in $C^\infty(\dot T^*
X,\mathrm{End}(\pi^*(E)))$. If the principal symbol is scalar and real
(\ref{transform}) is the transformation law of a partial connection
\footnote{A partial connection along a vector field $v$ can
be  defined by its covariant derivative which is a map $\nabla_v:
C^\infty(X;E) \to C^\infty(X;E)$ satisfying $\nabla_v(f g)=v(f)g+f \nabla_v g$
for all $f \in C^\infty(X), g \in C^\infty(X;E)$. Hence, parallel transport is defined
along $v$ only. Moreover, where $v$ vanishes this is a bundle homomorphism.}
along the vector field $H_P$. Hence, by
$\nabla_{H_P}:=H_P + \I \mathrm{sub}(P)$ a covariant derivative is defined.
We have therefore proved
the following proposition.

\begin{pro}
  Let $E$ be a vector bundle and suppose that
  $P \in \PsDO_{cl}^m(X,\Lambda^{n/2}X \otimes E)$ has real scalar principal
  symbol. Then the subprincipal symbol of $\mathrm{sub}(P)$
  defined locally by (\ref{subprince}) is invariantly defined as a partial
  connection on $\pi^* E$ along the Hamiltonian vector field $H_P$.
\end{pro}

Let $X$ be an oriented closed $n$-dimensional Riemannian manifold, let $E \to X$
be a hermitian vector bundle and suppose that
$P: C^\infty(X;E) \to C^\infty(X;E)$ is a formally selfadjoint
second order differential operator of Laplace type,
i.e. $\sigma_P(\xi)=g(\xi,\xi)\cdot \mathbf{1}$.
Then there exists a hermitian connection (see e.g. \cite{MR1215720} ch. 2.1)
\begin{gather*}
 \nabla: C^\infty(X;E) \to C^\infty(X;E \otimes T^*X)
\end{gather*}
 and a potential $V \in
 C^\infty(X;\mathrm{End}(E))$ such that
 \begin{gather*}
  P=\nabla^* \nabla + V
 \end{gather*}
The connection and the potential are uniquely determined by these properties.
The operator $P$ is essentially selfadjoint on $C^\infty(X;E)$
and in case it is positive we may define the square root
$P^{1/2}$ by functional calculus. By Seeley (\cite{MR0237943}) we know that
$P^{1/2}$ is a classical pseudodifferential operator
of order $1$ and its principal symbol is given by
$\sigma_{P^{1/2}}(\xi)=||\xi||_g \cdot \mathbf{1}$.
We use the metric to identify the bundle $\Lambda^{n/2}X$ with the trivial
bundle and in this way we understand $P$ and $P^{1/2}$ as operators
acting on $C^\infty(X;E \otimes \Lambda^{n/2}X)$.
Since $P$ is of Laplace type the Hamiltonian vector field
$H_{P^{1/2}}$ when restricted to the unit tangent bundle $T_1^*X$
coincides with the geodesic spray. We will therefore write
$H_g$ for $H_{P^{1/2}}$
 in order to emphasize the dependence from the metric.

Now $\nabla$ determines uniquely a hermitian partial connection
$\tilde \nabla_{H_g}$ along $H_g$ on $\pi^* E$ which satisfies
\begin{gather} \label{paral}
 \left(\tilde \nabla_{H_g} \pi^*(f)\right)(x,\xi) = \nabla_{\pi_*(H_g(x,\xi))} f \quad
 \textrm{for all}\; f \in C^\infty(X,E),
\end{gather}
where $\pi^*(f) \in C^\infty(T^*X;\pi^*(E))$ is the pull back of a section $f
\in C^\infty(X;E)$.
If we fix a local trivialization of $E$ we have
\begin{gather} \label{explicitconn}
   (\tilde \nabla_{H_g} f)(x,\xi) - (H_g f)(x,\xi)= \frac{\I}{||\xi||} \sum_{i,k}  g^{ik} A_i \xi_k,
\end{gather}
where $\nabla_{i}=\partial_i + \I A_i$.
The geometric meaning of this partial connection is as follows.
The vector field $H_g$ generates the geodesic flow. Hence, a partial connection
along $H_g$ allows to transport vectors along geodesics.
The partial connection (\ref{paral}) is chosen in such a way that
$v \in \pi^*(E)_{(x,\xi)}=E_x$ gets transported along the geodesic
with the original connection $\nabla$ on $E$.
We have
\begin{pro} \label{subdet}
 For $P=\nabla^*\nabla + V$ as above the partial connection determined by the
 subprincipal symbol $\mathrm{sub}(P^{1/2})$ of $P^{1/2}$
 coincides with $\tilde \nabla_{H_g}$.
\end{pro}
\begin{proof}
 We calculate everything in local coordinates where $|g|=1$ to
 keep the formulas as simple as possible.
 In such local coordinates we have $\nabla_i=\partial_i+ \I A_i$.
 Then one easily calculates
 \begin{gather*}
   \mathrm{sub}(P)(\xi)=\sum_{i,k} 2 \I g^{ik} A_i \, \xi_k.
 \end{gather*}
 Since the principal symbol of $P$ is scalar the formula as proved in
 \cite{MR0405514}
 \begin{gather}
   \mathrm{sub}(P^{1/2})(\xi)=\frac{1}{2} \sigma_P^{-\frac{1}{2}} \mathrm{sub}(P),
 \end{gather}
 continues to hold and
 we obtain
 \begin{gather*}
   \mathrm{sub}(P^{1/2})(\xi)= \I ||\xi||^{-1} \sum_{i,k} g^{ik} A_i \, \xi_k.
 \end{gather*}
 This coincides with the claimed formula.
\end{proof}

\subsection{Egorov's theorem}

Suppose that $A \in \PsDO_{cl}^1(X;E \otimes \Lambda^{n/2})$ has real scalar principal part,
let $H_A$ be the associated Hamiltonian vector field on $\dot T^* X$ and
let $\nabla_{H_A}=H_A + \I\, \mathrm{sub}(A)$ be the partial connection
on $\pi^*E$ defined by the subprincipal symbol.
Then this determines a geometric flow $\alpha_t$ on the vector bundle
$\pi^* E$ such that the flow lines $(x(t),\xi(t),v(t))$
consists of the orbits $(x(t),\xi(t))$
of the Hamiltonian flow and $v(t)$ expressed in coordinates of a local bundle chart satisfies
\begin{gather}
 \frac{dv(t)}{dt}=\I \mathrm{sub}(A) v(t).
\end{gather}
Note that $\alpha_t$ lifts the Hamiltonian
flow $h_t$ on $\dot T^* X$ and makes $\pi^* E$ an $\mathbb{R}$-equivariant vector bundle.
The induced flow $\alpha_t^*$
on $C^\infty(\dot T^* X,\pi^* E)$ satisfies
\begin{gather}
 \frac{d}{dt} \alpha_t^* f= \nabla_{H_A} f,
\end{gather}
which shows that the flow is defined independent of a choice of local coordinates.
Hence, there is also an action $\mathrm{Ad}(\alpha_t)$ of
$\mathbb{R}$ on $\pi^* \mathrm{End}(E)$ which extends the Hamiltonian
flow and is compatible with $\alpha_t$. This defines a flow on
$C^\infty(\dot T^* X,\pi^* \mathrm{End}(E))$ which we denote by
$\mathrm{Ad}(\alpha_t)^*$. Clearly, if $f \in C^\infty(\dot T^* X,\pi^*
\mathrm{End}(E))$
 \begin{gather}\label{bla1}
 \frac{d}{dt} \mathrm{Ad}(\alpha_t)^* f=  [ \nabla_{H_A}, f].
\end{gather}
In local coordinates one has
 \begin{gather} \label{bla2}
 [ \nabla_{H_A}, f]=H_A f + \I [\mathrm{sub}(A),f].
\end{gather}
Egorov's theorem now reads as follows.
\begin{pro}
 Let $A \in \PsDO_{cl}^1(X;E \otimes \Lambda^{n/2})$ and suppose that the principal
 symbol of $A$ is of real scalar type, i.e.
 $\sigma_A(\xi) = h(\xi) \cdot \mathbf{1}$, where $h \in C^\infty(\dot T^* M,\mathbb{R})$.
 Then, if $B$ is in $\PsDO_{cl}^m(X;E \otimes \Lambda^{n/2})$, also
 $B_t:=e^{+\I t A} B e^{-\I t A}$
 is in $\PsDO_{cl}^m(X;E)$.
 Moreover,
 \begin{gather*}
  \sigma_{B_t}=\mathrm{Ad}(\alpha_t)^*(\sigma_B).
 \end{gather*}
\end{pro}
\begin{proof}
 As usual we have
 \begin{gather} \label{flofe}
  \frac{d}{dt}B_t=\I [A,B].
 \end{gather}
The right hand side is a pseudodifferential operator in
$\PsDO_{cl}^n(X;E \otimes \Lambda^{n/2})$ and its
principal symbol is given by
\begin{gather} \label{flof}
  \sigma_{\I [A,B]}= \{\sigma_A,\sigma_B\}+ \I [\mathrm{sub}(A),B].
\end{gather}
This follows immediately from the formulas for the asymptotic
expansion of products of pseudodifferential operators. By
(\ref{bla1}) and (\ref{bla2}) equation (\ref{flofe}) is on the level
of principal symbols the flow equation for
$\mathrm{Ad}(\alpha_t)^*$. This equation can be solved on the symbol
level order by order and one can construct a classical symbol for
$B_t$ in the usual way (as for example carried out in
\cite{MR618463}, Ch VIII, \S 1). The only difference to the scalar
case is the second term in (\ref{flof}) which causes the Hamiltonian
flow to be replaced by $\mathrm{Ad}(\alpha_t)^*$.
\end{proof}
Suppose now that $X$ is an oriented Riemannian manifold and
let $A=P^{1/2}$, where $P=\nabla^* \nabla +V$ is a positive
Laplace type operator. We may use the metric to identify $\Lambda^{1/2}X$
with the trivial bundle. Moreover, the Hamiltonian
vector field $H_A$ restricted to the unit cotangent bundle $T_1^*X$
coincides with the geodesic spray. In this case it is convenient
to identify positively homogeneous functions on $\dot T^*X$ with
smooth functions on $T_1^*X$ by restriction. Hence, the
principal symbol of an operator in $\PsDO_{cl}^m(X;E \otimes \Lambda^{n/2})$
is in $C^\infty(T_1^*X,\pi^* \mathrm{End}(E))$.
In this case the Egorov theorem says, that
the flow $\alpha_t$ transports vectors in $\pi^* E$ parallel
with respect to the connection $\nabla$ along the geodesic flow on $T_1^*X$.

\section{The Dirac operator and the frame flow}

In this section $(X,g)$ is a compact oriented Riemannian manifold of dimension $n \geq 3$.
A spin structure (see e.g. \cite{MR1031992,MR1777332})
on $X$ is an $\Spin(n)$-principal bundle $P$ over $X$ together with a smooth
covering $\eta$ from $P$ onto the bundle $FX$ of oriented orthonormal frames,
such that the following diagram is commutative.
\begin{equation}
 \begin{CD}
   P \times \Spin(n) @>>> P @>>> X \\
   @VV \eta \times \lambda  V  @VV \eta V   @| \\
   FX \times \SO(n)   @>>> FX @>>> X
 \end{CD}
\end{equation}
Here $\lambda$ denotes the covering map $\Spin(n) \to \SO(n)$. The
(complexified) Clifford algebra $\mathrm{Cl}_c(\mathbb{R}^n)$ is
isomorphic to $\mathrm{Mat}(2^{[\frac{n}{2}]},\mathbb{C})$ if $n$ is
even and to $\mathrm{Mat}(2^{[\frac{n}{2}]},\mathbb{C}) \oplus
\mathrm{Mat}(2^{[\frac{n}{2}]},\mathbb{C})$ if $n$ is odd. The
Clifford modules $\Delta_n$ are then defined by the action of this
matrix algebra on $\mathbb{C}^{2^{[\frac{n}{2}]}}$, where in case
$n$ is odd we project onto the first summand. Hence, $\Delta_n$ is
an irreducible module for the Clifford algebra
$\mathrm{Cl}_c(\mathbb{R}^n)$. Since $\Spin(n) \subset
\mathrm{Cl}_c(\mathbb{R}^n)$ the Clifford modules $\Delta_n$ are
also modules for the group $\Spin(n)$. The corresponding
representation $\rho: \Spin(n) \to \mathrm{Aut}(\Delta_n)$ is called
the spinor representation of $\Spin(n)$. This representation is
irreducible if $n$ is odd. It is the direct sum of two irreducible
components if $n$ is even. The spinor bundle $S$ associated with a
Spin structure is the associated bundle $P \times_\rho \Delta_{n}$.
The Levi-Civita connection on $FX$ lifts naturally to a connection
on $P$ and this defines a connection $\nabla_S: C^\infty(X;S) \to
C^\infty(X;S \otimes T^* X)$ on $S$, the Levi-Civita connection on
the spinor bundle. The Dirac operator $D: C^\infty(X;S) \to
C^\infty(X;S)$ is defined by $-\I \gamma \circ \nabla_S$, where
$\gamma$ denotes the action of covector fields on sections of the
spinor bundle by Clifford multiplication. $D$ is of Dirac type and
essentially selfadjoint on $C^\infty(X;S)$. The operator
$F=\mathrm{sign}(D)$ is in $\PsDO_{cl}^0(X;S)$ and its principal
symbol $\sigma_F(\xi)$ is given by Clifford multiplication by
$\frac{1}{|\xi|}\xi$. If $n$ is even, then Clifford multiplication
with the volume form times $\I^{\frac{n(n+1)}{2}}$ defines an
involution $\Gamma$ on $L^2(X;S)$  which anti-commutes with $D$ and
with $F$, but which commutes with $|D|$. The Lichnerowicz formula
allows us to express the square of the Dirac operator by the spinor
Laplacian
\begin{gather}
 D^2=\nabla^*_S \nabla_S +\frac{1}{4}R,
\end{gather}
where $R$ is the scalar curvature.
By Prop. \ref{subdet} the connection determined by the subprincipal
symbol $\mathrm{sub}(|D|)$ of $|D|$ transports spinors along geodesics
with the spinor connection $\nabla_S$. The corresponding flow
on the bundle $\pi^* S \to T^*_1 X$ will be denoted by $\alpha_t$. The induced flow
on the bundle $\pi^* \mathrm{End}(S) \to T^*_1 X$ will be denoted as before by
$\mathrm{Ad}(\alpha_t)$.
This induces a $1$-parameter group $\beta_t=\mathrm{Ad}(\alpha_t)^*$  of
$*$-automorphisms of the algebra $C^\infty(T_1^*X,\pi^* \mathrm{End}(S))$. It
extends continuously to the $C^*$-algebra $\mathcal{A}=C(T_1^*X,\pi^* \mathrm{End}(S))$.
The Egorov theorem of the previous section therefore reads as follows.
\begin{pro}
 Let $D$ be the Dirac operator on a compact spin manifold $X$ and let $A \in \PsDO_{cl}^0(X;S)$.
 Then with $A_t:=e^{+\I t |D|} A e^{-\I t |D|}$ we have $A_t \in
 \PsDO_{cl}^0(X;S)$ for all $t \in \rbb$ and
 \begin{gather*}
   \sigma_{A_t}=\beta_t(\sigma_A).
 \end{gather*}
\end{pro}
\noindent
 And as a consequence we get
\begin{cor}\label{invstate}
 Let $X$ be a closed Riemannian spin manifold of dimension $n \geq 3$
 and let $D$ be the Dirac operator.
 Let $\phi_k$ be a sequence of eigensections to $D^2$ with
 \begin{gather*}
  D^2 \phi_k = \lambda_k \phi_k,\quad \langle \phi_k,\phi_j \rangle=\delta_{kj}.
 \end{gather*}
 such that $|\lambda_k| \nearrow \infty$ and such that
 the sequence of states $$\omega_i(A):=\langle \phi_i, A \phi_i \rangle$$
 on the $C^*$-algebra $\overline{\PsDO_{cl}^0(X;S)}$
 converges in the weak-$*$-topology.
 Then there is a $\beta_t$-invariant state $\omega_\infty$ on
 $C(T_1^* X,\pi^* \mathrm{End}(S))$ such that
 \begin{gather*}
  \lim_{n \to \infty} \omega_n(A)=\omega_\infty(\sigma_A)
 \end{gather*}
\end{cor}
\begin{proof}
 Since the states $\omega_i$ are invariant under the flow induced by conjugation
 with $e^{\I |D| t}$ so is the limit state.
 The limit state $\omega$ vanishes on all operators of order $-1$ since their
 product with $|D|$ is bounded.
 Since the norm closure of those operators are the compact operators
 $\mathcal{K}$ we have $\omega(\mathcal{K})=\{0\}$.
 Hence, the state projects to a state on the quotient
 $\overline{\PsDO_{cl}^0(X;S)}/\mathcal{K}$. The symbol map is known
 to extend by continuity to a map
 $\overline{\PsDO_{cl}^0(X;S)} \to C(T_1^* X,\pi^* \mathrm{End}(S))$
 with kernel $\mathcal{K}$. Hence, the quotient is via the symbol
 map isomorphic to $C(T_1^* X,\pi^* \mathrm{End}(S))$. Moreover, by the above the
 symbol map is equivariant with respect to the two flows.
\end{proof}

There are two natural invariant states $\omega_{\pm}$ on $C(T_1^*
X,\pi^* \mathrm{End}(S))$ given by
\begin{gather}
\omega_\pm(a):= \frac{2}{\mathrm{rk}(S)}\int_{T_1^* X}
\mathrm{Tr}(P_{\pm}(\xi) a(\xi) P_{\pm}(\xi)) d \mu_L(\xi),
\end{gather}
where $\mu_L$ is the normalized Liouville measure and
$P_{\pm}:=\frac{1}{2}(1 \pm \sigma_F)$. The tracial state
\begin{gather*}
 \omega=\frac{1}{2}(\omega_+ + \omega_-),\\
 \omega(a)=\frac{1}{\mathrm{rk}(S)} \int_{T_1^* X} \mathrm{Tr}(a(\xi)) d \mu_L(\xi)
\end{gather*}
is therefore \underline{not} ergodic (see Appendix \ref{appa} for the notion of ergodicity in this context).
\begin{theorem}\label{ifergodic}
 Suppose the frame flow is ergodic. Then $\omega_+$ and $\omega_-$ are ergodic
 states with respect to $\beta_t$. If moreover $n$ is odd the systems
 $\left(C(T_1^* X,\pi^* \mathrm{End}(S)),\omega_\pm\right)$ are $\beta_t$-abelian.
\end{theorem}
\begin{proof}
 Note that since the Levi-Civita connection is compatible with
 the Clifford multiplication both $P_+$ and $P_-$ are easily seen to be
 invariant under $\beta_t$.
 To prove ergodicity we have to show (see Appendix \ref{appa})
 that all $\beta_t$-invariant elements in
 $P_\pm L^\infty(T^*_1 X,\pi^* \mathrm{End}(S)) P_\pm$ are of the form $c P_\pm$
 with $c \in \mathbb{C}$. To show $\beta_t$-abelianness we have to show
 that an invariant element in $L^2(T^*_1 X,\pi^* \mathrm{End}(S))$
 is of the form $c_1 P_+ + c_2 P_-$ with $c_1,c_2 \in \mathbb{C}$.
 We therefore analyze, how invariant elements in $L^2(T^*_1 X,\pi^* \mathrm{End}(S))$
 look like. \\
 {\bf Step 1:}
 Denote by $\mathbb{C}\langle X_1,\ldots,X_{n} \rangle$ the space of noncommutative
 polynomials in the variables $X_1,\ldots,X_{n}$.
 Now we define a continuous map
 \begin{gather}
   T: L^2(FX) \otimes \mathbb{C}\langle X_1,\ldots,X_{n} \rangle
   \to L^2(T^*_1 X,\pi^* \mathrm{End}(S)),\\
   T(f \otimes p)(\xi):=\int_{F_{\xi}X} f(\xi,v) p(\xi,v) d\mu(v),
 \end{gather}
 where the integration is over the invariant measure on the fiber $F_{\xi}X$.
 The action of covectors on spinors here is by Clifford multiplication.
 The pullback $h^*_t$ of the frame flow on $FX$ defines a flow
 on $L^2(FX) \otimes \mathbb{C}\langle X_1,\ldots,X_{n} \rangle$
 by acting on the first tensor factor.
 Since the connection is compatible with Clifford multiplication and the Clifford
 action the map $T$ intertwines the pullback $h^*_t$ of the frame flow and the
 flow $\beta_t$, i.e.
 \begin{gather}
  T \circ (h_t^* \otimes \mathbf{1})=\beta_t \circ T.
 \end{gather}
 Note that the Clifford action on the spinor bundle is irreducible. This implies
 that $T$ is onto.\\
 {\bf Step 2:}
 Let $\Psi \in L^2(T^*_1 X,\pi^* \mathrm{End}(S))$ be an invariant vector.
 By the above we can choose an element
 $f \in L^2(FX) \otimes \mathbb{C}\langle X_1,\ldots,X_{n} \rangle$ such that
 $T(f) = \Psi$. Since by assumption the frame flow is
 ergodic we have $\lim_{T \to \infty} \frac{1}{T}\int_{0}^T h^*_t(f) dt =
 \int_{FX} f(x) dx= m$, where $m$ is a constant polynomial.
 Since $\Psi$ is invariant and $T$ is equivariant
 we have almost everywhere
 \begin{gather*}
  \Psi(\xi)= T(m)(\xi)=\left( \int_{F_\xi X} m(\xi,v) d\mu(v) \right).
 \end{gather*}
 Since the measure on $F_\xi X$ is invariant under the action of $SO(n-1)$
 the endomorphism $\int_{F_\xi X} m(\xi,v) d\mu(v)$ commutes with the action of
 $\Spin(n-1)$. Clearly, the projections $P_+(\xi)$ and $P_-(\xi)$
 commute with the $\Spin(n-1)$ action, and it is easy to see (for example by
 calculating the dimensions) that the representations of $\Spin(n-1)$
 on the ranges of $P_+(\xi)$ and $P_-(\xi)$ are irreducible. Moreover,
 these two representations are equivalent iff $n$ is even. In this case
 the algebra of invariant matrices in $\mathrm{End(S_\xi)}$ is generated by
 $P_+(\xi),P_-(\xi),\Gamma$. This shows that an invariant element
 of the form $P_\pm A P_\pm$ is proportional to $P_\pm$ which proves ergodicity.
 If $n$ is odd the two irreducible representations of $\Spin(n-1)$
 on the ranges of $P_+(\xi)$ and $P_-(\xi)$ are inequivalent. Hence, any invariant
 element in $A \in \mathrm{End(S_\xi)}$ is of the form $c_1 P_+ + c_2 P_-$.
\end{proof}
\noindent
{\bf Remark:}
 That $\left(C(T_1^* X,\pi^* \mathrm{End}(S)),\omega_\pm\right)$
 is not $\beta_t$-abelian in even dimensions has a simple interpretation. It is that the space
 of invariant states on $C(T_1^* X,\pi^* \mathrm{End}(S))$ is not a simplex, i.e.
 the decomposition of an invariant state into ergodic states is not unique.
 For example the tracial state $\omega$ has the decomposition
 $\omega=\frac{1}{2}(\omega_+ + \omega_-)$, but it also has the decomposition
 $\omega=\frac{1}{2}(\omega_1 + \omega_2)$ where $\omega_1(a)=\omega((1+\Gamma) a)$
 and $\omega_2(a)=\omega((1-\Gamma) a)$. The states $\omega_1$ and $\omega_2$
 can also be shown to be ergodic if the frame flow is ergodic.

\begin{theorem}\label{erg:spin}
 Let $X$ be a compact Riemannian spin manifold with spinor bundle $S$ and Dirac
 operator $D: C^\infty(X;S) \to C^\infty(X;S)$.
 Let $\phi_k$ be an orthonormal sequence of eigensections
 of $D$ with eigenvalues $\lambda_k \nearrow \infty$ such that $\phi_k$
 spans\footnote{in the sense that the hull of the vectors is dense.}
 $L^2_+(X;S)=\frac{1+F}{2}L^2(X;S)$.
 If the frame flow of $X$ is ergodic, then Quantum Ergodicity
 holds, i.e.
 \begin{gather*}
  \lim_{N \to \infty} \frac{1}{N}\sum_{k=1}^N |\langle \phi_k, A \phi_k
  \rangle-\omega_+(\sigma_A)|
  = 0, \quad \textrm{for all}\; A \in \PsDO_{cl}^0(X;S).
 \end{gather*}
 Moreover, there is a density one
 subsequence $\phi_k'$ such that
 \begin{gather*}
  \lim_{k \to \infty} \langle \phi_k', A \phi_k' \rangle
  = \omega_+(\sigma(A)), \quad \textrm{for all}\; A \in \PsDO_{cl}^0(X;S).
 \end{gather*}
 An analogous statement holds for eigensections with $\lambda_k \searrow -\infty$
 with $\omega_+$ replaced by $\omega_-$.
\end{theorem}
\begin{proof}
We denote $\omega=\frac{1}{2}(\omega_+ +
\omega_-)=\frac{1}{\mathrm{rk}(S)}\int_{T^*_1 X} \tr(\sigma_A(\xi))
d\mu(\xi)$ be the ordinary tracial state on $C^\infty(T^*_1 X,\pi^*
\mathrm{End}(S))$. The heat trace asymptotics (cf. \cite{GrSee})
\begin{gather}
\tr(A \cdot e^{-D^2 t}) = C(n) \int_{T^*_1 X} \tr(\sigma_A(\xi))
d\mu(\xi) \cdot t^{-n/2} + O(t^{-n/2+1/2})
\end{gather}
 which one easily gets from the classical calculus of pseudodifferential
 operators together with Karamata's Tauberian theorem implies that
 \begin{gather} \label{themean}
\lim_{N \to \infty} \frac{1}{N}\sum_{k=1}^N \langle \frac{1+F}{2}
\phi_k, A \frac{1+F}{2}\phi_k \rangle
  = \omega_+(\sigma_A), \qquad A \in \PsDO_{cl}^0(X;S).
 \end{gather}
 The proof is now an analog to the proof of Shnirelman, Zelditch and Colin
 de Verdi\`ere (\cite{MR0402834, MR1239173, MR0818831, MR0916129}).
 To prove Quantum ergodicity it is obviously enough to show that
 for any selfadjoint $A \in \PsDO_{cl}^0(X;S)$ with $\omega_+(\sigma_A)=0$ we have
 \begin{gather}
  \lim_{N \to \infty }\frac{1}{N}\sum_{k=1}^N |\langle \phi_k, A \phi_k \rangle|=0.
 \end{gather}
 Clearly, $|\langle \phi_k, A \phi_k \rangle|=
 |\langle \phi_k, \frac{1+F}{2} A \frac{1+F}{2} \phi_k \rangle|$ and therefore, we assume
 without loss of generality that $A=\frac{1+F}{2} A \frac{1+F}{2}$. Hence,
 $\sigma_A=P_+ \sigma_A P_+$.
 Now let $A_T:=\frac{1}{T} \int_{0}^T e^{+\I t |D|} A e^{-\I t |D|} dt$
 and since $\omega_+$ is ergodic we have
 $\omega_+(|\sigma_{A_T}|^2) \to 0$ (see Lemma \ref{thelemma}) as $T \to \infty$. This implies
 that also $\omega_+(|\sigma_{A_T}|) \to 0$ as $T \to \infty$.
 On the other hand we have
 \begin{gather}
  |\langle \phi_k, A \phi_k \rangle|=|\langle \phi_k, A_T \phi_k \rangle| \leq
  \langle \phi_k, |A_T| \phi_k \rangle.
 \end{gather}
 Then formula (\ref{themean}) applied to $|A_T|$ together with the
 fact that $\omega_+(\sigma_{|A_T|}) < \epsilon$ for $T$ large enough
 allows us to conclude that
 \begin{gather}
  \limsup_{N \to \infty }\frac{1}{N}\sum_{k=1}^N |\langle \phi_k, A \phi_k
  \rangle| < \epsilon,
 \end{gather}
 for all $\epsilon>0$ from which the assertion follows.
 The existence of a density one sequence is based on a diagonalization
 argument which is well known (see e.g. \cite{MR1384146}).
\end{proof}

\section{The Hodge-Laplace}\label{sec:hodge}

Let $X$ be a compact oriented Riemannian manifold, let $\Lambda^* X:=\Lambda^* T^* X$
be the exterior algebra bundle, and let $\Lambda^*_\cbb X$ be its complexification.
We denote by $d$ the exterior derivative,
$\delta$ the coderivative (i.e. the formal adjoint of $d$), and by $*$ the Hodge
star operator. Then $d+\delta$ is a Dirac type operator
acting on sections of $\Lambda^*_\cbb X$. Its square is the Hodge Laplace
operator
\begin{gather}
 \Delta = d \delta + \delta d.
\end{gather}
Note that $\Delta$ leaves the subspace of $p$-forms invariant. In the
following we will denote the restriction of $\Delta$ to $p$-forms
by $\Delta_p$. The ordinary Laplace operator acting on functions
is therefore equal to $\Delta_0$.
If $\nabla_p$ is the Levi-Civita covariant derivative of $p$-forms
the Weitzenb\"ock formula states that
\begin{gather}
 \Delta_p = \nabla_p^* \nabla_p + H_p,
\end{gather}
where $H_p$ is section of  $\mathrm{End}(\Lambda^p T^* X)$ which can be
expressed in terms of the curvature of the connection.
For example $H_1$ is equal to the Ricci curvature.
We conclude that the partial connection determined by the subprincipal
symbol of $\Delta^{1/2}_p$ transports $p$-multivectors along the Hamiltonian vector
field parallel with respect to the Levi-Civita connection. The corresponding
flow on the $C^*$-algebra $C(T_1^*X,\pi^* \mathrm{End}(\Lambda^p_\cbb X))$
will be denoted by $\beta_t$. There is a natural invariant tracial state $\omega_{\mathrm{tr}}$
on $C(T_1^*X,\pi^* \mathrm{End}(\Lambda^p_\cbb X))$ given by
\begin{gather}
 a \to \omega_{\mathrm{tr}}(a)= { n \choose p }^{-1} \int_{T_1^* X} \tr(a(\xi)) dL(\xi).
\end{gather}
As in the Dirac case this state is not ergodic for $0< p <n$.
Let $P \in C(T_1^*X,\pi^* \mathrm{End}(\Lambda^p_\cbb X))$ defined
by
\begin{gather}
 P(\xi) v:= i(\xi) \xi \wedge v,
\end{gather}
where $i(\xi)$ denotes the operator of interior multiplication with $\xi$.
Then $P$ is an orthogonal projection
in $C(T_1^*X,\pi^* \mathrm{End}(\Lambda^p_\cbb X))$ which is invariant
under $\beta_t$, and hence
\begin{gather}
 \omega_{\mathrm{tr}}=\frac{n-p}{n} \,\omega_t+ \frac{p}{n}\, \omega_l,\\
 \omega_l(a)=\frac{n}{p}\,\omega_{\mathrm{tr}}((1-P) \cdot a),\\
 \omega_t(a)=\frac{n}{n-p} \,\omega_{\mathrm{tr}}(P \cdot a),
\end{gather}
is a decomposition into invariant states.
The non-ergodicity of this state can be seen as the classical counterpart
of the Hodge decomposition
\begin{gather}
 C^\infty(X;\Lambda^p_\cbb X)=d  C^\infty(X;\Lambda^{p-1}_\cbb X) \oplus
 \delta C^\infty(X;\Lambda^{p+1}_\cbb X) \oplus \mathrm{ker}(\Delta_p),
\end{gather}
which induces a decomposition of $L^2(X;\Lambda^p_\cbb X)$ into invariant
subspaces of $\Delta_p$.

\begin{theorem}\label{p:erg}
 Suppose $0<p<n, \; n \geq 3$ and let
 $k_p=2\, \mathrm{min}(p,n-p)$. Suppose that the $k_p$-frame flow is ergodic
 and that $p \not= \frac{n-1}{2}$.
 Then $\omega_t$ is an ergodic
 state with respect to $\beta_t$ on $C(T_1^* X,\pi^* \Lambda^p_\cbb X)$. If moreover
 $p \not= \frac{n}{2}$ then the system
 $\left(C(T_1^* X,\pi^* \Lambda^p_\cbb X),\omega_t \right)$ is $\beta_t$-abelian.
\end{theorem}

\begin{proof}
 The proof is similar to the proof in the Dirac case.
 We investigate the set of invariant vectors in
 $P L^2(T_1^* X,\pi^* \mathrm{End}(\Lambda^p_\cbb X)) P$
 and in $L^2(T_1^* X,\pi^* \mathrm{End}(\Lambda^p_\cbb X)) P$.\\
 {\bf Step 1:}
 Let $N=k_p$. Denote by $\mathcal{P} \subset \mathbb{C}
 \langle X_1,\ldots,X_N, Y_1,\ldots,Y_N\rangle$
  the ring of polynomials in the noncommutative variables
  $X_1,\ldots,X_N$ and $Y_1,\ldots,Y_N$ such that in each summand the same number of $X$
  and $Y$ occur.
  Now define a map
  \begin{gather*}
   T: L^2(F_{N}X) \otimes \mathcal{P}
   \to L^2(T_1^*X,\pi^* \mathrm{End}(\Lambda^p_\cbb X)),\\
   T(f \otimes p)(\xi) = \int_{F_{N,\xi} X} f(\xi,v) \hat p(v) d\mu(v),
  \end{gather*}
  where integration is over the fibre $F_{N,\xi} X$ of the bundle $F_{N} X$
  over the point $\xi \in T_1^*X$. The endomorphism $\hat p(v)$ is defined by
  replacing all $X_i$ by exterior multiplication with $v_i$ and all
  $Y_i$ by interior multiplication with $v_i$. Since the number of $X$ and $Y$
  is the same in each summand the operators leave the space
  $\Lambda^p_\cbb T^*_{\pi(\xi)} X$ invariant.
  Since exterior and interior multiplication are compatible with the Levi-Civita
  connection the map $T$ intertwines the pullback of $N$-frame flow and the flow
  $\beta_t$, i.e.
  \begin{gather*}
   T \circ (h_t^* \otimes \mathbf{1})=\beta_t \circ T.
  \end{gather*}
  Moreover, by an elementary exercise in linear algebra any endomorphism of
  $\Lambda^p_\cbb T^*_{\pi(\xi)} X$ can be represented by a linear combination
  of elements of the form $\hat p(v_1,\ldots,v_{N})$, where $v_1,\ldots,v_{N}$
  is a frame and $p$ has degree at most $2N$. Therefore, the map $T$ is surjective.\\
 {\bf Step 2:}
   Now let $\Psi$ be an invariant element in
  $L^2(T_1^* X,\pi^* \mathrm{End}(\Lambda^p_\cbb X))$. Then we may find an
  $f \in L^2(F_{N}X) \otimes \mathcal{P}$ such that $\Psi=T(f)$.
  By the same argument as in the proof of theorem \ref{ifergodic}
  it follows that $\Psi=T(m)$, where $m$ is some constant polynomial
  such that $\hat m(\xi)$ commutes with the action of $SO(n-1)$
  on each fiber of $\pi^*(\Lambda^p_\cbb X)$ . Note that $P(\xi)$ and $(1-P(\xi))$ project
  onto invariant subspaces.
  The fiber of $\pi^*(\Lambda^p_\cbb X)$ at the point $(x,\xi) \in T^*_1 X$
  is given by $\Lambda^p_\cbb T^*_x X$. On the other hand $T^*_x X = \rbb \xi \oplus V$, where
  $V$ is the orthogonal complement of $\xi$ in $T^*_x X$. Hence,
  we have the decomposition $\Lambda^p T^*_x X= \mathbb{R} \otimes
  \Lambda^{p-1} V \oplus \Lambda^{p} V$.
  It is now easy to see that $1-P$ projects onto $\mathbb{C} \otimes \Lambda^{p-1}_\cbb V$, whereas
  $P$ projects onto $\Lambda^{p}_\cbb V$. The representation of $SO(n-1)$
  on $\Lambda^{p} \mathbb{C}^{n-1}$ is irreducible (see e.g. \cite{MR1153249}, Lecture 18)
  since we assumed $p\not=\frac{n-1}{2}$. If $p \not=\frac{n}{2}$ then
  all other components which occur in the decomposition into irreducible
  representations are inequivalent to this representation as one
  can see by calculating the dimensions (the exceptional case $n=3, p=2$
  by other methods). Hence, in this case the algebra
  of invariant elements in $\mathrm{End}(\Lambda^p_\cbb T^*_x X)$ is generated
  by $P(\xi)$ and $1-P(\xi)$. Therefore, any invariant element is of the
  form $c_1 P + c_2 (1-P)$. This shows $\mathbb{R}$-abelianness and ergodicity if
  $p \not= \frac{n}{2}$.
  Ergodicity for $p=\frac{n}{2}$ follows from the fact that any invariant
  element in $P(\xi) \mathrm{End}(\Lambda^p_\cbb T^*_x X) P(\xi)$ is proportional
  to $P(\xi)$ which is a simple consequence of the irreducibility of
  the $SO(n-1)$ action on $\Lambda^p_\cbb V$.
\end{proof}

The orthoprojection onto the closure of
$\delta  C^\infty(X;\Lambda^{p-1} T^* X)$ is given by the zero order
pseudodifferential operator $(\Delta_p|_{\mathrm{ker}(\Delta_p)^{\perp}})^{-1}\delta d$.
The principal symbol of this operator is easily seen to coincide with $P$.
Because of the Hodge decomposition quantum ergodicity for the operator
$\Delta_p$ with $0<p<n$ can never hold in the strict sense.
An eigenform $\phi$ of $\Delta_p$ with nonzero eigenvalue can be decomposed
uniquely as
\begin{gather}
 \phi= d \phi_- + \delta \phi_+,
\end{gather}
where $\phi_-$ is an eigenform of $\Delta_{p-1}$ and
$\phi_+$ is an eigenform of $\Delta_{p+1}$. Hence, part of the spectrum
of $\Delta_p$ comes from part of the spectrum of $\Delta_{p-1}$. The other part
can be obtained by solving the system
\begin{gather} \label{system}
 \Delta_p \phi = \lambda \phi,\\ \nonumber
 \delta \phi =0.
\end{gather}
We will show that in certain situations this system is quantum ergodic.
\begin{theorem}\label{p:qerg}
 Assume $0<p<n,\; n\geq 3$ and
 suppose that $\phi_k$ is an orthonormal sequence of eigen-$p$-forms satisfying
 \begin{gather*}
 \Delta_p \phi_k = \lambda_k \phi_k,\\
 \delta \phi_k =0,
\end{gather*}
such that the $\phi_k$ span $\mathrm{ker} (\delta)$ and $\lambda_k \nearrow \infty$.
Suppose that $p \not= \frac{n-1}{2}$. Then, if the
$2\,\mathrm{min}(p,n-p)$-frame flow is ergodic,
the system is quantum ergodic in the sense that
\begin{gather*}
  \lim_{N \to \infty }\frac{1}{N} \sum_{k \leq N} |\langle \phi_k, A \phi_k \rangle -
 \omega_t(\sigma_A)| =0,
 \end{gather*}
 for all $A \in \PsDO_{cl}^0(X;\Lambda^p_\cbb X)$. In particular there is a density one subsequence
 $\phi_k'$ such that
 \begin{gather*}
  \lim_{k \to \infty} \langle \phi_k', A \phi_k' \rangle
  = \omega_t(\sigma_A), \quad \textrm{for all}\; A \in \PsDO_{cl}^0(X;\Lambda^p_\cbb X).
 \end{gather*}
\end{theorem}
\begin{proof}
 The proof is along the same lines as the proof of theorem \ref{erg:spin}.
 Suppose that $A \in \PsDO_{cl}^0(X;\Lambda^p_\cbb X)$ with $\omega_t(A)=0$ and $A^*=A$.
 Let $F$ be the operator $(\Delta_p|_{\mathrm{ker}(\Delta_p)^{\perp}})^{-1} \delta d$,
 then we have $\sigma_F=P$.
 Since $F \phi_k = \phi_k$ we conclude that
 $\langle \phi_k, A \phi_k \rangle=\langle \phi_k, F A F \phi_k \rangle$.
 Hence, we may assume without loss of generality
 that $A=FAF$ and hence, $\sigma_A=P \sigma_A P$. Again we have the heat asymptotics
 \begin{gather}
  \tr (A e^{-\Delta_p t}) \sim C(p,n) \omega_{\mathrm{tr}}(\sigma_A) t^{-n/2},
 \end{gather}
 and from Lemma \ref{thelemma} $\omega_t(|\sigma_{A_T}|^2) \to 0$ as $T \to \infty$.
 Together these statements with the Karamata's Tauberian theorem imply Quantum ergodicity
 in the stated form exactly in the same way as in the proof of theorem \ref{erg:spin}.
\end{proof}

In the above proof it was necessary to exclude the case
$p=\frac{n-1}{2}$ because in this case the representation of
$SO(n-1)$ on $\Lambda^p \mathbb{C}^{n-1}$ is not irreducible but
splits into a direct sum of two irreducible representations. One
reason for this is the existence of an involution defined by the
Hodge star operator which commutes with the $SO(n-1)$ action. This
actually causes the state $\omega_t$ to be non-ergodic in case
$p=\frac{n-1}{2}$. To see this let $P_\pm \in C(T_1^*X,\pi^*
\mathrm{End}(\Lambda^p X))$ defined by
\begin{gather}
 P_\pm(\xi) v:= \frac{1}{2}\left(1 \pm \I^{p} i(\xi) *\right) i(\xi) \xi \wedge v,
\end{gather}
where $*: \Lambda^{p-1}_\cbb X \to \Lambda^{p+1}_\cbb X$ is the Hodge star operator.
Then $P_\pm$ are orthogonal projections
$C(T_1^*X,\pi^* \mathrm{End}(\Lambda^p_\cbb X))$ which commute in addition with $P$.
Therefore, we have $P=P_+ + P_-$ and
\begin{gather*}
 \omega_t = \frac{1}{2}(\omega_+ + \omega_-),\\
 \omega_+(a)=2 \omega_t(P_+ a),\\
 \omega_+(a)=2 \omega_t(P_- a).
\end{gather*}
Again these states are invariant. The same proof as for theorem \ref{p:erg}
now gives
\begin{theorem}
 Suppose that $n>1$ is odd. Let $p=\frac{n-1}{2}$. If the $(n-1)$-frame flow is ergodic
 then $\omega_\pm$ are ergodic
 states with respect to $\beta_t$ on $C(T_1^* X,\pi^* \Lambda^p_\cbb X)$ and
 the systems $\left(C(T_1^* X,\pi^* \Lambda^p_\cbb X),\omega_\pm \right)$ are $\beta_t$-abelian.
\end{theorem}

It is amusing that also in case $p=\frac{n-1}{2}$ the non-ergodicity
of the state $\omega_t$ is related to the existence of a ``quantum
symmetry'', i.e. of a pseudodifferential operator which commutes
with $\Delta_p$ and leaves the kernel of
$\overline{\mathrm{Rg}(\delta)}$ invariant. Namely, for
$p=\frac{n-1}{2}$  the operator $\I^{p+1} \Delta_p^{-1/2} \delta *$
is a selfadjoint involution on $\overline{\mathrm{Rg}(\delta)}$
whose principal symbol is precisely $\I^p i(\xi) *$. Hence, $P_\pm$
are the principal symbols of the projections to the $\pm 1$
eigenspaces of this involution. The proof of theorem \ref{p:qerg}
gives
\begin{theorem}
 Suppose that $n$ is odd and $2p=n-1$.
 Suppose that $\phi_k$ is an orthonormal sequence of eigen-$p$-forms satisfying
 \begin{gather*}
 \Delta_p \phi_k = \lambda_k \phi_k,\\
 \delta \phi_k =0,\\
 \I^{p+1} \delta * \phi_k = \pm \sqrt{\lambda_k} \phi_k
\end{gather*}
such that the $\phi_k$ span $\overline{\mathrm{Ran}(\delta \pm
\I^{p+1}\Delta_p^{-1/2}\delta * \delta)}$ and with $\lambda_k
\nearrow \infty$. Then, if the $(n-1)$-frame flow is ergodic, the
system is quantum ergodic in the sense that
\begin{gather*}
\lim_{N \to \infty }\frac{1}{N} \sum_{k \leq N} |\langle \phi_k, A
\phi_k \rangle - \omega_\pm(\sigma_A(\xi))| =0,
\end{gather*}
for all $A \in \PsDO_{cl}^0(X;\Lambda^p_\cbb X)$. In particular
there is a density one subsequence $\phi_k'$ such that
\begin{gather*}
  \lim_{k \to \infty} \langle \phi_k', A \phi_k' \rangle
  = \omega_\pm(\sigma(A)), \quad \textrm{for all}\; A \in \PsDO_{cl}^0(X;\Lambda^p_\cbb X).
\end{gather*}
\end{theorem}

\appendix
\section{Ergodicity of states for noncommutative classical systems} \label{appa}

Let $\mathcal{A}$ be a unital $C^*$-algebra.
A state $\omega$ over $\mathcal{A}$ is a positive linear functional
$\omega: \mathcal{A} \to \mathbb{C}$ with $\omega(1)=1$. The set of states $E_\mathcal{A}$
of $\mathcal{A}$ is a convex weakly-$*$-compact subset of $\mathcal{A}^*$
and its extreme points are the pure states $P_\mathcal{A}$.
Each state $\omega$ gives via the GNS-construction rise to a representation
$\pi_\omega: \mathcal{A} \to \mathcal{L}(\mathcal{H}_\omega)$
with a cyclic vector $\Omega_\omega$ such that
$\omega(a)=\langle \Omega_\omega,\pi_\omega(a) \Omega_\omega \rangle$.
Up to equivalence the triple $(\pi_\omega,\mathcal{H}_\omega,\Omega_\omega)$
is uniquely determined by its properties and we refer to it as the GNS triple.

In the following let $\alpha_t$
be a strongly continuous group of $*$-auto\-morphisms of $\mathcal{A}$.
The set of invariant states $E_\mathcal{A}^{\alpha_t}$ is again a weakly-$*$-compact
subset of $\mathcal{A}^*$. The extreme points in $E_\mathcal{A}^{\alpha_t}$ are called
ergodic states. Hence, an invariant state $\omega$ is ergodic if it cannot
be written as a convex linear combination of two other invariant states.
If $\omega$ is an invariant state the group $\alpha_t$ can be uniquely implemented
by a strongly continuous unitary group $U(t)$ on the GNS-Hilbert space
$\mathcal{H}_\omega$ such that $U(t) \Omega_\omega=\Omega_\omega$
and $U(t)^* \pi_\omega(a) U(t) = \pi_\omega(\alpha_t(a))$.
Let $E_\omega$ be the orthogonal projection onto the space of $\alpha_t$-invariant
vectors in $\mathcal{H}_\omega$. Then the pair $(\mathcal{A},\omega)$ is
called $\mathbb{R}$-abelian if all operators in $E_\omega \pi_\omega(\mathcal{A}) E_\omega$
commute pairwise. If we look at the following conditions
\begin{enumerate}
 \item $E_\omega$ has rank one,
 \item $\omega$ is ergodic for $\alpha_t$, i.e. $\omega \in
 \mathcal{E}(E_\mathcal{A}^{\alpha_t})$,
 \item $\{\pi_\omega(\mathcal{A}) \cup U(t)\}$ is irreducible on $\mathcal{H}_\omega$,
\end{enumerate}
then it is known
that $(1) \Rightarrow (2) \Leftrightarrow (3)$.
If, moreover, either $(\mathcal{A},\omega)$ is $\mathbb{R}$-abelian, or
$\Omega_\omega$ is separating for $\pi_\omega(\mathcal{A})''$, then
all three conditions are equivalent (see \cite{MR545651}, Prop. 4.3.7, Th. 4.3.17, Th 4.3.20).

Now let $E \to X$ be a hermitian complex vector bundle over a
compact Hausdorff space $X$. Then $\mathcal{A}=C(X;\mathrm{End}(E))$
is a unital $C^*$-algebra. Suppose that $\mu$ is some finite Borel
measure with $\mu(X)=1$ and let $P \in C(X;\mathrm{End}(E))$ be a
non-trivial orthogonal projection onto a subbundle, i.e. $P$ has
constant rank $k$. Then $\omega(a):=\frac{1}{k} \int_X \tr\left(P(x)
a(x) P(x)\right) d\mu(x)$ is a state over $\mathcal{A}$. The GNS
triple can be calculated explicitly and is given by
$\mathcal{\mathcal{H}}_\omega=L^2(X;\mathrm{End}(E)) P$ with scalar
product $\langle a,b \rangle = \frac{1}{k}\int_X \tr(a^*(x) b(x))
d\mu(x)$ and $\Omega_\omega=P$. The action of $\mathcal{A}$ on
$\mathcal{\mathcal{H}}_\omega$ is by multiplication from the left.
The von Neumann closure of $\pi_\omega(\mathcal{A})$ is given by
$\pi_\omega(\mathcal{A})''=L^\infty(X;\mathrm{End}(E))$. The
commutant $\pi_\omega(\mathcal{A})'$ of $\pi_\omega(\mathcal{A})$
can be identified with the opposite algebra of $P
L^\infty(X;\mathrm{End}(E)) P$ which acts on $L^2(X;\mathrm{End}(E))
P$ by right multiplication. Note that $\Omega_\omega$ is separating
for $\pi_\omega(\mathcal{A})''$ iff $P=\mathrm{Id}$. Now any
continuous geometric flow on $E$ determines a continuous geometric
flow on $\mathrm{End}(E)$. If the hermitian structure is preserved
by the flow, this gives rise to a strongly continuous $1$-parameter
group $\alpha_t$ on $\mathcal{A}$. If $P$ and $\mu$ are invariant
under the flow, then $\omega$ is an invariant state. If all
invariant vectors in $L^2(X;\mathrm{End}(E)) P$ are of the form $c
P$ with $c \in \mathbb{C}$, then by the above $\omega$ is ergodic.
In this case $E_\omega \pi_\omega(\mathcal{A}) E_\omega$ is just a
multiplication by a number and the system is $\mathbb{R}$-abelian.
If a state $\omega'$ is majorized by $\omega$ then there exists a
positive element $A \in P L^\infty(X;\mathrm{End}(E))P$ such that
such that $\omega'(a)=\omega(a A)$. Hence, if $P$ is up to a
constant the only invariant element in
$PL^\infty(X;\mathrm{End}(E))P$ then $\omega$ is ergodic.

Now denote by $\omega'$ the restriction of the state $\omega$
to the subalgebra $\mathcal{B}=P C(X;\mathrm{End}(E))P$.
If $E'$ denotes the subbundle onto which $P$ projects then clearly
$\mathcal{B}=C(X;\mathrm{End}(E'))$ and $\omega'$
becomes the tracial state on $\mathcal{B}$. Ergodicity of the state is equivalent
to the condition that all invariant elements in
$P L^\infty(X;\mathrm{End}(E))P=L^\infty(X;\mathrm{End}(E'))$
are multiples of $P$. Therefore, ergodicity of $\omega'$
is equivalent to the ergodicity of $\omega$. Since
$\Omega_{\omega'}$ is separating for $\pi''_{\omega'}(\mathcal{B})$
ergodicity of $\omega$ is equivalent to $E_{\omega'}$
having rank one. As a consequence we get
\begin{lem}\label{thelemma}
 If $A \in PL^\infty(X,E)P$ such that $\omega(A)=0$, then
 ergodicity of $\omega$ implies that
 \begin{gather}
  \lim_{T \to \infty} \omega(|A_T|^2) = 0,
 \end{gather}
 where
 \begin{gather*}
  A_T=\frac{1}{T} \int_0^T \alpha_t(A) dt.
 \end{gather*}
\end{lem}
\begin{proof}
 Since $A = PAP$ we also have $A_T^2=P A_T^2 P$. By the above $E_{\omega'}$
 has rank one and its range is spanned by $P$. Now by the von Neumann
 ergodic theorem $A_T$ converges to $E_{\omega'} A$ in
 $P L^2(X,E) P$. But since $\omega(A)=0$ we have $E_{\omega'} A=0$
 and therefore $A_T$ converges to $0$ in $P L^2(X,E) P$.
 By definition this means that $\omega(|A_T|^2)$ converges to 0.
\end{proof}

\noindent{\bf Acknowledgements.} The authors would like to thank J.
Bolte for bringing the problem to their attention and for many
useful discussions.  The authors would also like to thank M. Brin,
D. Dolgopyat, P. Gerard, N. Kamran, M. Lesch, W. M\"uller, M. Pollicott, I.
Polterovich, P. Sarnak, A. Shnirelman, R. Schubert, J. Stix, J. Toth
and S. Zelditch for useful discussions.  The auhtors would also like 
to thank the anonymous referee whose comments helped improve the 
presentation of the results.  

The first author would like
to thank the Department of Theoretical Physics at the university of
Ulm for their hospitality during his visit. This paper was started
while the first author visited Max Planck Institute for Mathematics
in Bonn, and completed while he was visiting IHES; their
hospitality is greatly appreciated. The second
author would like to thank McGill university and CRM Analysis
laboratory for the hospitality during his Montreal visit.



\bigskip

\begin{thebibliography}{BGV92}

\bibitem[Arn61]{MR0158330}
V. I. Arnold.
\newblock Some remarks on flows of line elements and frames.
\newblock {\em Dokl. Akad. Nauk SSSR} {\bf 138} (1961), 255--257.

\bibitem[BGV92]{MR1215720}
N. Berline, E. Getzler, and M. Vergne.
\newblock {\em Heat kernels and {D}irac operators}, volume {\bf 298} of {\em
Grundlehren der Mathematischen Wissenschaften [Fundamental
Principles of Mathematical Sciences]}.
\newblock Springer-Verlag, Berlin, 1992.

\bibitem[BoK98]{MR1644088} J. Bolte and S. Keppeler.
\newblock Semiclassical time evolution and trace formula for relativistic
spin-$1/2$ particles.
\newblock {\em Phys. Rev. Lett.} {\bf 81} (1998), no. 10, 1987--1991.

\bibitem[BoK99]{MR1694732}  J. Bolte and S. Keppeler.
\newblock A semiclassical approach to the Dirac equation.
\newblock {\em Ann. Physics} {\bf 274} (1999), no. 1, 125--162.

\bibitem[Bol01]{MR1838911} J. Bolte.
\newblock Semiclassical expectation values for relativistic
particles with spin 1/2. Invited papers dedicated to Martin C. Gutzwiller, Part III.
\newblock {\em Found. Phys.} {\bf 31} (2001), no. 2, 423--444.

\bibitem[BoG04]{MR2073612} J. Bolte and R. Glaser.
\newblock Zitterbewegung and semiclassical
observables for the Dirac equation.
\newblock {\em J. Phys. A} {\bf 37} (2004), no. 24, 6359--6373.

\bibitem[BoG04.2]{BoG04.2} J. Bolte, R. Glaser.
\newblock A semiclassical Egorov theorem
and quantum ergodicity for matrix valued operators.
\newblock {\em Comm. Math. Phys.} 247 (2004), no. 2, 391--419.

\bibitem[BR79]{MR545651}
O. Bratteli and D.~W. Robinson.
\newblock {\em Operator algebras and quantum statistical mechanics. {V}ol. 1}.
\newblock Springer-Verlag, New York, 1979.
\newblock $C\sp{\ast} $- and $W\sp{\ast} $-algebras, algebras, symmetry groups,
  decomposition of states, Texts and Monographs in Physics.

\bibitem[Br75]{MR0370660}
M. Brin.
\newblock Topological transitivity of one class of dynamical systems
and flows of frames on manifolds of negative curvature.
\newblock {\em Funct. Anal. Appl.} {\bf 9}, 8--16 (1975)

\bibitem[Br76]{MR0394764}
M. Brin.
\newblock The topology of group extensions of Anosov systems.
\newblock {\em Math. Notes} {\bf 18}, 858--864, (1976)

\bibitem[Br82]{MR0670078}
M. Brin.
\newblock Ergodic theory of frame flows.
\newblock In {\em Ergodic Theory and Dynamical Systems II,} Proc.
Spec. Year, Maryland 1979-80, Progr. Math. {\bf 21}, 163--183,
Birkh\"auser, Boston, 1982.

\bibitem[BrG80]{MR0582702}
M. Brin and M. Gromov.
\newblock On the ergodicity of frame flows.
\newblock {\em Inv. Math.} {\bf 60}, 1--7 (1980)

\bibitem[BrK84]{MR0756723}
M. Brin and H. Karcher.
\newblock Frame flows on manifolds with pinched negative curvature.
\newblock {\em Comp. Math.} {\bf 52}, 275--297 (1984)

\bibitem[BrP74]{MR0343316} M.~I. Brin and  Ja.~B. Pesin.
{\em Partially hyperbolic dynamical systems.}
 Izv. Akad. Nauk SSSR Ser. Mat. {\bf 38} (1974), 170--212.

\bibitem[BuP03]{MR1988429}
K. Burns and M. Pollicott.
\newblock Stable ergodicity and frame flows.
\newblock {\em Geom. Dedicata}, {\bf 98}, 189--210 (2003)

\bibitem[CV85]{MR0818831}
Y. Colin de Verdi\`ere.
\newblock  Ergodicit\'e et fonctions propres du laplacien.
\newblock {\em  Comm. Math. Phys.}  {\bf 102}, 497--502, (1985)

\bibitem[D82]{MR661876}
N. Dencker.
\newblock  On the propagation of polarization sets for systems of real
   principal type.
\newblock {\em J. Funct. Anal.}  46:351--372, 1982.

\bibitem[DG75]{MR0405514}
J.~J. Duistermaat and V.~W. Guillemin.
\newblock The spectrum of positive elliptic operators and periodic
  bicharacteristics.
\newblock {\em Invent. Math.}, {\bf 29(1)}, 39--79, (1975)

\bibitem[DH72]{Hormander:1972}
J.~J. Duistermaat and L.~H{\"o}rmander.
\newblock Fourier integral operators. {II}.
\newblock {\em Acta Mathematica}, {\bf 128}, 183--269, (1972)

\bibitem[EW96]{EW96} G. Emmrich and A. Weinstein.
\newblock Geometry of the transport
equation in multicomponent WKB approximations.
\newblock {\em Comm. Math. Phys.} 176 (1996), no. 3, 701--711.

\bibitem[Fri]{MR1777332}
T. Friedrich.
\newblock {\em Dirac operators in {R}iemannian geometry}, volume~25 of {\em
  Graduate Studies in Mathematics}.
\newblock American Mathematical Society.
\newblock Translated from the 1997 German original by Andreas Nestke.

\bibitem[FuHa91]{MR1153249}
W. Fulton and J. Harris.
\newblock {\em Representation theory}, volume~129 of {\em
  Graduate Texts in Mathematics}.
\newblock Springer Verlag, New York, 1991.

\bibitem[GrSee95]{GrSee}
G. Grubb and R. Seeley.
\newblock Weakly parametric pseudodifferential operators and
Atiyah-Patodi-Singer boundary problems.
\newblock {\em Invent. Math.} {\bf 121}, 481--529, (1995)

\bibitem[GMMP97]{GMMP97}
P. G\'erard, P. Markowich, N. Mauser and F. Poupaud.
\newblock Homogenization limits and Wigner transforms.
\newblock {\em Comm. Pure Appl. Math.} 50 (1997), no. 4, 323--379.
Erratum:  {\em Comm. Pure Appl. Math.} 53 (2000), no. 2, 280--281.

\bibitem[LM89]{MR1031992}
H.~B. Lawson, Jr. and M. Michelsohn.
\newblock {\em Spin geometry}, volume~38 of {\em Princeton Mathematical
  Series}.
\newblock Princeton University Press, Princeton, NJ, 1989.

\bibitem[San99]{San99}
M.~R. Sandoval.
\newblock Wave-trace asymptotics for operators of Dirac type.
\newblock {\em Comm. PDE} 24  (1999), no. 9-10, 1903--1944.

\bibitem[See67]{MR0237943}
R.~T. Seeley.
\newblock Complex powers of an elliptic operator.
\newblock In {\em Singular Integrals (Proc. Sympos. Pure Math., Chicago, Ill.,
  1966)}, pages 288--307. Amer. Math. Soc., Providence, R.I., 1967.

\bibitem[Shn74]{MR0402834}
A.~I. Shnirelman.
\newblock Ergodic properties of eigenfunctions.  (Russian) .
\newblock {\em  Uspehi Mat. Nauk}  {\bf 29}, 181--182, (1974).

\bibitem[Shn93]{MR1239173} A.~I. Shnirelman.
\newblock On the asymptotic properties of eigenfunctions in the
regions of chaotic motion.
\newblock In V. Lazutkin {\em KAM theory and semiclassical
approximations to eigenfunctions.} Ergebnisse der Mathematik und
ihrer Grenzgebiete (3), 24. Springer-Verlag, Berlin, 1993.


\bibitem[Tay81]{MR618463}
M.~E. Taylor.
\newblock {\em Pseudodifferential operators}, volume~34 of {\em Princeton
  Mathematical Series}.
\newblock Princeton University Press, Princeton, N.J., 1981.

\bibitem[Zel87]{MR0916129}
S. Zelditch.
\newblock Uniform distribution of eigenfunctions on compact hyperbolic surfaces.
\newblock {\em  Duke Math. J.} {\bf 55}, 919--941 (1987)

\bibitem[Zel96]{MR1384146}
S. Zelditch.
\newblock Quantum ergodicity of $C^*$ dynamical systems.
\newblock {\em Comm. Math. Phys.} {\bf 177}, 507--528, (1996)



\end{thebibliography}
\end{document}